\documentclass[a4paper,10pt]{scrartcl}

\usepackage[left=3.5cm, right=3.5cm, top=2.3cm, bottom=2.1cm]{geometry}
\usepackage[utf8]{inputenc}
\usepackage{babel}
\usepackage[T1]{fontenc}
\usepackage{amsmath}
\usepackage{amssymb}
\usepackage{multicol}
\usepackage{graphicx}
\usepackage{xcolor}
\usepackage{lmodern}
\usepackage{float}
\usepackage{caption}
\DeclareMathSizes{10}{10}{8.3}{7}

\title{On an acyclic relaxation of incomparable families of sets}
\author{Maximilian Krone}
\date{October 9, 2024}

\newcounter{ThNr}[section]
\DeclareRobustCommand{\newTh}[1]{\refstepcounter{ThNr}\arabic{ThNr} \label{#1}}

\begin{document}
\parskip3explus1exminus1ex
\parindent0mm

\begin{center}
\textbf{\Large On an acyclic relaxation of incomparable families of sets}\\ \vspace*{5mm}

{\large Maximilian Krone}

Technische Universität Ilmenau\\
June 17, 2025 \vspace*{2mm}
\end{center}

\begin{abstract}
For two families $\mathcal{A}, \mathcal{B} \subseteq \mathcal{P}([k])$, we write $\mathcal{A}\vdash\mathcal{B}$ if $A\not\supseteq B$ for each two sets $A \in \mathcal{A}$ and $B \in \mathcal{B}$. $\mathcal{A}$ and $\mathcal{B}$ are called \textit{incomparable} if $\mathcal{A}\vdash\mathcal{B}$ and $\mathcal{B}\vdash\mathcal{A}$. Seymour proved that the maximum size of two incomparable equal-sized families in $\mathcal{P}([k])$ is $\frac{1}{4}2^k$.

A sequence of families $\mathcal{B}_1,\dots,\mathcal{B}_l \ \subseteq \mathcal{P}([k])$ is called $d$-\textit{exceeding} if $\mathcal{B}_i\vdash\mathcal{B}_j$ for all $i,j\in [l]$ with $j-i\in [d]$. Cyclically reusing $d+1$ pairwise incomparable families yields arbitrarily long $d$-exceeding sequences of families. We prove inversely that the maximum size of equal-sized families of a sufficiently long $1$-exceeding sequence in $\mathcal{P}([k])$ is also $\frac{1}{4}2^k$. 

A sequence of sets $B_1,\dots,B_l \subseteq [k]$ is called $d$-\textit{exceeding} if $\{B_1\},\dots,\{B_l\}$ is $d$-exceeding, that is, if $B_i \not\supseteq B_j$ for all $i,j\in [l]$ with $j-i\in [d]$. We locate the maximum $d$ such that there exist arbitrarily long $d$-exceeding sequences of subsets of $[k]$ between $(1-o(1)) \tfrac{1}{e} 2^k$ and $\tfrac{1}{2}2^k-2$.
\end{abstract}

\hrulefill

\section*{Introduction}

Two sets $A,B\subseteq [k]=\{1,\dots,k\}$ are called \textit{incomparable} if $A\not\supseteq B$ and $B\not\supseteq A$. Two families of sets $\mathcal{A}, \mathcal{B} \subseteq \mathcal{P}([k])$ are called \textit{incomparable} if each two sets $A \in \mathcal{A}$, $B \in \mathcal{B}$ are incomparable. 
In \cite{Seymour}, Seymour proved that the maximum size of two equal-sized incomparable families in $\mathcal{P}([k])$ is $\frac{1}{4}2^k$. In \cite{false}, Gerbner, Lemons, Palmer, Patkós and Szécsia presented a simple generalization of Seymour's construction to a larger number of pairwise incomparable families and conjectured it to be optimal. Recently, in \cite{New}, Behague, Kuperus, Morrison and Wright disproved the conjecture by presenting a better construction.

In Section \ref{incomparable families}, we present a combination of both constructions focusing on families of the same size. This yields pairwise incomparable families whose number can be polynomial in $k$ or even growing faster (like in \cite{false}) but whose size is larger (like in \cite{New}). In particular, for $k=\omega(1)$ and $1\leq a = o\left(\sqrt{\frac{k}{\log_2 k}}\right)$, there exist $r \sim \frac{1}{a!}\left(\frac{k}{\log_2 k}\right)^a$ pairwise incomparable families in $\mathcal{P}([k])$ of size $\sim \frac{a^a}{e^a a!} \frac{2^k}{r}$ (see Corollary \ref{incomparable families}.\ref{double asymptotic}). The symbol $\sim$ states equality except for a factor $(1\pm o(1))$. With $a=1$, we achieve the construction from \cite{New}.

The main goal of the present paper is to introduce an acyclic relaxation of the concept of pairwise incomparable families, which avoids cyclic conditions $A_1\not\supseteq A_2 \not\supseteq \dots \not\supseteq A_s \not\supseteq A_1$.
For two families $\mathcal{A}, \mathcal{B} \subseteq \mathcal{P}([k])$, we write $\mathcal{A}\vdash\mathcal{B}$ if $A\not\supseteq B$ (that is, $B\setminus A \neq \emptyset$) for each two sets $A \in \mathcal{A}$ and $B \in \mathcal{B}$. Hence, $\mathcal{A}$ and $\mathcal{B}$ are incomparable if and only if $\mathcal{A}\vdash\mathcal{B}$ and $\mathcal{B}\vdash\mathcal{A}$.
A sequence of families $\mathcal{B}_1,\dots,\mathcal{B}_l \ \subseteq \mathcal{P}([k])$ is called $d$-\textit{exceeding} if $\mathcal{B}_i\vdash\mathcal{B}_j$ for all $i,j\in [l]$ with $j-i\in [d]$. For $d=1$, the sequence is called \textit{exceeding}.

Cyclically reusing $r$ pairwise incomparable families yields arbitrarily long $(r-1)$-exceeding sequences of families. In Section \ref{exceeding families}, we also prove some inverse results: For fixed $r$, some upper bounds on the size of $r$ pairwise incomparable equal-sized families are also true for sufficiently long $(r-1)$-exceeding sequences. In particular (for $r=2$), the maximum size of the families in $\mathcal{P}([k])$ in a sufficiently long exceeding sequence is $\frac{1}{4}2^k$, which is equal to Seymour's bound from \cite{Seymour} on the size of two incomparable families.\newpage

In Section \ref{exceeding sets}, we study the special case of families of size $1$:
A sequence of sets $B_1,\dots,B_l \subseteq [k]$ is called $d$-\textit{exceeding}, if $\{B_1\},\dots,\{B_l\}$ is $d$-exceeding, that is, if $B_i \not\supseteq B_j$ for all $i,j\in [l]$ with $j-i\in [d]$.
For $k\geq 1$, let $\delta(k)$ be the maximum $d$ such that, for all $l\in \mathbb{N}$, there exists a $d$-exceeding sequence of length $l$ of subsets of $[k]$. Clearly, $\delta(1)=0$. We determine the first values $\delta(2)=1$, $\delta(3)=2$ and $\delta(4)=5$, and prove the general asymptotic bounds
$$(1-o(1))\tfrac{1}{e}2^k\ \leq\ \delta(k)\ \leq\ \tfrac{1}{2}2^k -2\ .$$
In \cite{self}, these result are used to examine the existence of cut covers of powers of directed paths.

\section{Pairwise incomparable families of sets} \label{incomparable families}

In this Section, we present a very general construction of pairwise incomparable families of the same size. The idea for the construction has already been introduced in \cite{New}. We start with a simple idea from \cite{false}.

\textbf{Lemma \ref{incomparable families}.\newTh{doubling lemma}\ -\ Upscaling pairwise incomparable families \cite{false}\\}
Let $1\leq s\leq k$. If there are $r$ pairwise incomparable families in $\mathcal{P}([s])$ of size $\beta 2^s$, then there are $r$ pairwise incomparable families in $\mathcal{P}([k])$ of size $\beta 2^k$. \enlargethispage{5mm}

\textbf{Proof. \ } Let $\mathcal{B}_1,\dots,\mathcal{B}_r \subseteq \mathcal{P}([s])$ of size $\beta 2^s$ be pairwise incomparable. Define $\mathcal{B}'_1,\dots,\mathcal{B}'_r \subseteq \mathcal{P}([k])$ as
$$\mathcal{B}'_i := \big\{B \cup A\ |\ B \in \mathcal{B}_i, A\subseteq [k]\setminus [s] \big\}.$$
Clearly, $|\mathcal{B}'_i|=\beta 2^s 2^{k-s} = \beta 2^k$. Let $B_i \in \mathcal{B}'_i$ and $B_j \in \mathcal{B}'_j$ for $i\neq j$. Then, $B_i\cap [s] \in \mathcal{B}_i$ and $B_j\cap [s] \in \mathcal{B}_j$, so they are incomparable. Hence, also $B_i$ and $B_j$ are incomparable \cite{false}. \hfill $\Box$\\

In \cite{false}, this observation is used for the following simple construction of incomparable families in $\mathcal{P}([k])$:
Let $0 \leq s \leq k$. By starting with the $\binom{s}{\lfloor s/2 \rfloor}$ pairwise incomparable families each containing exactly one unique subsets of $[s]$ of size $\lfloor s/2 \rfloor$, it follows from Lemma \ref{incomparable families}.\ref{doubling lemma} that there exist $\binom{s}{\lfloor s/2 \rfloor}$ many pairwise incomparable families in $\mathcal{P}([k])$ of size $\frac{1}{2^s}2^k$. This construction is optimal in the extreme cases $s=k$ (Sperner's Theorem \cite{Sperner}) and $s=0$ (but clearly not for $s=1$). For $s=2$, it yields the two incomparable families $\big\{ B\subseteq [k] : B\cap \{1,2\} = \{1\} \big\}$ and $\big\{ B\subseteq [k] : B\cap \{1,2\} = \{2\} \big\}$ of size $\frac{1}{4} 2^k$. It was first proved by Seymour \cite{Seymour} that this size is best possible.

\begin{figure}[H]
\begin{center}
	\vspace*{-1mm}
	\includegraphics[height=3cm]{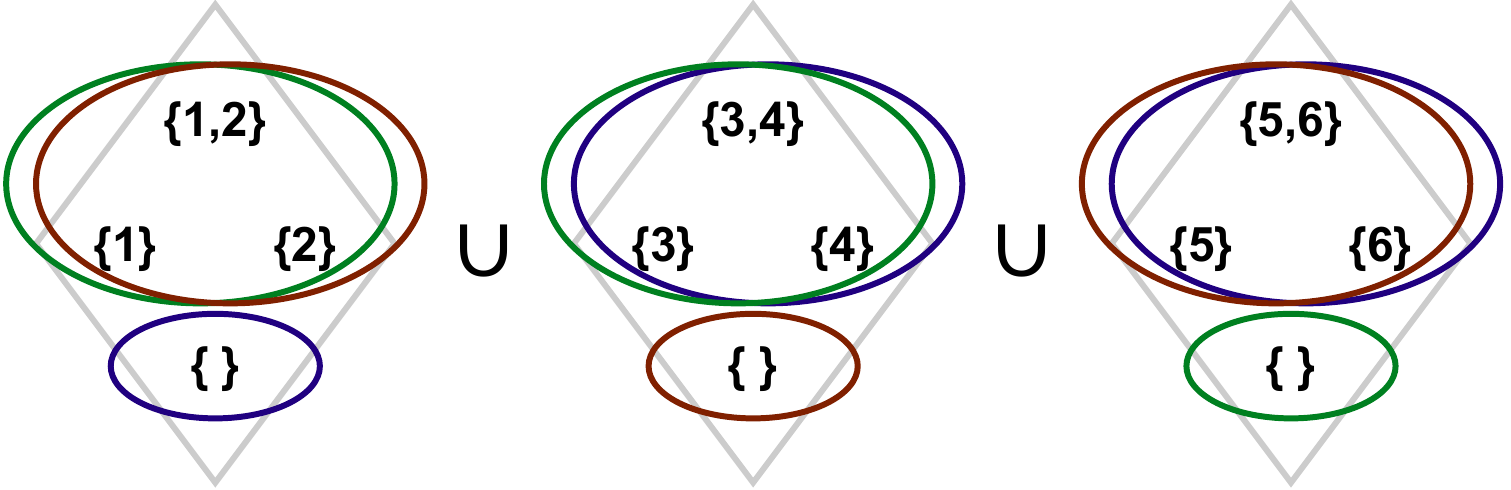}
	\captionsetup{justification=centering}
	\caption{Three pairwise incomparable families in $\mathcal{P}([6])$ of size $9$:\\ \vspace*{1mm}
		\textcolor[HTML]{008020}{$\big\lbrace\ A\cup B\cup C\ :\ A\in \big\lbrace \{1\}, \{2\}, \{1,2\} \big\rbrace,\quad B \in \big\lbrace \{3\}, \{4\}, \{3,4\} \big\rbrace,\quad C=\emptyset\ \big\rbrace$}\\ \vspace*{1mm}
		\textcolor[HTML]{802000}{$\big\lbrace\ A\cup B\cup C\ :\ A\in \big\lbrace \{1\}, \{2\}, \{1,2\} \big\rbrace,\quad B=\emptyset,\quad C\in \big\lbrace \{5\}, \{6\}, \{5,6\} \big\rbrace\ \big\rbrace$}\\ \vspace*{1mm}
		\textcolor[HTML]{200080}{$\big\lbrace\ A\cup B\cup C\ :\ A=\emptyset,\quad B \in \big\lbrace \{3\}, \{4\}, \{3,4\} \big\rbrace,\quad C\in \big\lbrace \{5\}, \{6\}, \{5,6\} \big\rbrace\ \big\rbrace$}}
	\label{3in6}
\end{center}
\end{figure}

For $s=3$, the simple construction yields three pairwise incomparable families of size $\frac{1}{8} 2^k$ in $\mathcal{P}([k])$. For $k\geq 6$, we can improve the fraction $\frac{1}{8}=0.125$ to $\frac{9}{64}\approx 0.141$, as shown in Figure \ref{3in6}. For $k=\omega(1)$, it will even be improved to $\frac{4}{27}\approx 0.148$ in Corollary \ref{incomparable families}.\ref{constant number}.
 \vspace{-8mm} \newpage

The simple construction gets even worse for a larger number of families. For $s=\omega(1)$, the proportion of the union of all families in $\mathcal{P}([k])$ is $\binom{s}{\lfloor s/2 \rfloor}\tfrac{1}{2^s} \sim \sqrt{\tfrac{2}{\pi s}} = o(1)$ (using Stirling's approximation). We will see in Corollary \ref{incomparable families}.\ref{double asymptotic} that we can achieve a constant proportion as long as the number of families is polynomial in $k$.

To reduce the number of parameters, we do not give true generalizations of the previous constructions, but versions, where some parameters are already optimized.

\textbf{Theorem \ref{incomparable families}.\newTh{exact}\ -\ Exact version\\}
Let $1 \leq a \leq \frac{s}{2}$ such that $\frac{s}{a}$ is a power of $2$, and $k \geq s \log_2 \frac{s}{a}$. Then, there exist $\binom{s}{a}$ pairwise incomparable families in $\mathcal{P}([k])$ of size $b=\frac{a^a(s-a)^{s-a}}{s^s} 2^k$.

\textbf{Proof. \ } Let $c=\log_2\frac{s}{a}$, that is, $2^c=\frac{s}{a}$. With usage of Lemma \ref{incomparable families}.\ref{doubling lemma}, we may assume that $k=s \log_2 \frac{s}{a} = sc$. In this case,
$$b = \tfrac{a^a(s-a)^{s-a}}{s^s}\ 2^{sc} = \tfrac{a^a(s-a)^{s-a}}{s^s} \left(\tfrac{s}{a}\right)^s = \left(\tfrac{s-a}{a}\right)^{s-a} = \left(2^c-1\right)^{s-a}.$$
We consider a partition of $[k]$ into $s$ sets $C_1,\dots,C_s$ of size $c$.
For each of the $\binom{s}{a}$ many sets $A\subseteq[s]$ of size $a$, define a family
$$\mathcal{B}_A = \left\{B \subseteq [k]\ |\ B\cap C_i = \emptyset\ \text{ if and only if } \ i\in A\right\}.$$
Indeed, $\mathcal{B}_A$ contains $(2^c-1)^{s-a}$ many sets, since for each $i \in [s]\setminus A$, we have $2^c-1$ possibilities to include or exclude the elements of $C_i$.

Finally, let $B \in \mathcal{B}_A$ and $B' \in \mathcal{B}_{A'}$, for distinct subsets $A\neq A'$ of $[s]$ of size $a$. We need to show that $B\not\supseteq B'$. There exists some element $i \in A\setminus A'$. From $i\in A$, we get $B\cap C_i = \emptyset$. From $i\notin A'$, we get $B'\cap C_i \neq \emptyset$. Hence $B\cap C_i \not\supseteq B'\cap C_i$, and $B\not\supseteq B'$. \hfill $\Box$\\

For the case that $\frac{s}{a}$ is not a power of $2$, we now prove an asymptotic version. Unfortunately, this demands a higher lower bound on $k$ for large $a$: For example, consider $a=\Theta(s)$. In the exact version, it suffices if $k$ is linear in $s$. In the asymptotic version, we need an extra factor $\log_2 s$ (which is indeed necessary for the construction, for example if $s=2a+1$).

\textbf{Theorem \ref{incomparable families}.\newTh{asymptotic}\ -\ Asymptotic version\\}
Let $1 \leq a \leq \frac{s}{2}$ and $s \big(\log_2 s - \tfrac{1}{2} \log_2 a +\omega(1) \big) \leq k$.
Then, there exist $\binom{s}{a}$ pairwise incomparable families in $\mathcal{P}([k])$ of size $b\sim\frac{a^a(s-a)^{s-a}}{s^s} 2^k$.

\textbf{Proof. \ } Let $c=\log_2 s - \tfrac{1}{2} \log_2 a +\omega(1) \in \mathbb{N}$, that is, $2^c=\omega(1)\frac{s}{\sqrt{a}}$, and let $d=\frac{a}{s}2^c$. By Lemma \ref{incomparable families}.\ref{doubling lemma}, we may assume that $k=sc$. 

Again, we consider a partition of $[k]$ into $s$ sets $C_1,\dots,C_s$ of size $c$. For each $i \in [s]$, let $\mathcal{D}_i \subseteq \mathcal{P}(C_i)$ be a family of $\lceil d \rceil$ sets of smallest possible size. That is, the sets in $\mathcal{D}_i$ have at most as many elements as the sets in $\mathcal{P}(C_i) \setminus \mathcal{D}_i$. This guarantees $D\not\supseteq D'$ for each two $D\in \mathcal{D}_i$ and $D'\in \mathcal{P}(C_i) \setminus \mathcal{D}_i$.
For each of the $\binom{s}{a}$ many sets $A\subseteq[s]$ of size $a$, define a family
$$\mathcal{B}_A = \left\{B \subseteq [k]\ |\ B\cap C_i \in \mathcal{D}_i\ \text{ if and only if } \ i\in A\right\}.$$
First, we check that the families are indeed pairwise incomparable. Let $B \in \mathcal{B}_A$ and $B' \in \mathcal{B}_{A'}$, for distinct subsets $A\neq A'$ of $[s]$ of size $a$. We need to show that $B\not\supseteq B'$. There exists some element $i \in A\setminus A'$. From $i\in A$, we get $B\cap C_i \in \mathcal{D}_i$. From $i\notin A'$, we get $B'\cap C_i \notin \mathcal{D}_i$. Hence $B\cap C_i \not\supseteq B'\cap C_i$, and $B\not\supseteq B'$.\\

Finally, we have to compute the size $b=|\mathcal{B}_A|$ of the families. Clearly, $b = \lceil d \rceil^a (2^c - \lceil d \rceil)^{s-a}$: For each $i\in A$, we have $\lceil d \rceil$ possibilities to include or exclude the elements of $C_i$. For each $i\in [s]\setminus A$, we have $2^c-\lceil d \rceil$ possibilities.

We calculate the value $b'$ of $b$, when we replace $\lceil d \rceil = d + \Delta$, for some $\Delta\in [0,1)$, by $d=\frac{a}{s}2^c$:
$$b' = d^a (2^c - d)^{s-a} = \left(\tfrac{a}{s}\right)^a \left(1-\tfrac{a}{s}\right)^{s-a} 2^{cs} = \tfrac{a^a(s-a)^{s-a}}{s^s} 2^k.$$
Let $z=\tfrac{\Delta s}{2^c}$. Hence, we have the conformity $\frac{z}{\Delta} = \frac{s}{2^c} = \frac{a}{d}$, so 
$$\tfrac{\lceil d \rceil}{d} = \tfrac{d+\Delta}{d} = \tfrac{a+z}{a},\ \text{ and }\ 
\tfrac{2^c-\lceil d \rceil}{2^c-d} = \tfrac{2^c-d-\Delta}{2^c-d} = \tfrac{s-a-z}{s-a}.$$
We now show that $b\sim b'$. The upper bound $b\leq b'$ follows from
$$\tfrac{b}{b'} = \left(\tfrac{\lceil d \rceil}{d}\right)^a \left(\tfrac{2^c-\lceil d \rceil}{2^c-d}\right)^{s-a} = \left(1+\tfrac{z}{a}\right)^a \left(1-\tfrac{z}{s-a}\right)^{s-a} \leq e^z e^{-z} = 1.$$
The other way around,
$$\tfrac{b'}{b} = \left(\tfrac{a}{a+z}\right)^a \left(\tfrac{s-a}{s-a-z}\right)^{s-a} = \left(1-\tfrac{z}{a+z}\right)^a \left(1+\tfrac{z}{s-a-z}\right)^{s-a},$$
which we can upper-bound by $e$ to the power of
$$ z \left(\tfrac{-a}{a+z} + \tfrac{s-a}{s-a-z}\right) = \tfrac{z^2s}{(a+z)(s-a-z)} \leq \tfrac{z^2s}{a(s/2)} = \tfrac{2z^2}{a}.$$
For the last inequality, we have used that $a+z\leq \tfrac{s}{2}$. This is a consequence of $a\leq \tfrac{s}{2}$, which translates to $d \leq \tfrac{1}{2}2^c \in \mathbb{N}$, and hence $d+\Delta = \lceil d \rceil \leq \tfrac{1}{2}2^c$. 

We obtain the (quite explicit) lower bound $b \geq e^{-2z^2/a}\ b' \geq \left(1-\tfrac{2z^2}{a} \right) b'$. 

For $z=\tfrac{\Delta s}{2^c}\leq \tfrac{s}{2^c} = o(\sqrt{a})$, this indeed yields $b \sim b' = \tfrac{a^a(s-a)^{s-a}}{s^s} 2^k$.  \hfill $\Box$\\

For $a=1$ and constant $s=r$, we conclude the following result. 

\textbf{Corollary \ref{incomparable families}.\newTh{constant number}\ -\ A constant number of families \cite{New}\\}
Let $r\geq 1$ be fixed. For $k=\omega(1)$, there exist $r$ pairwise incomparable families in $\mathcal{P}([k])$ of size $\sim\tfrac{(r-1)^{r-1}}{r^r} 2^k$. The fraction lies in the interval $\Big(\tfrac{1}{er},\,\tfrac{1}{e(r-1)}\Big)$:
$$\tfrac{1}{er} < \tfrac{1}{r}\left(1+\tfrac{1}{r-1}\right)^{-(r-1)} = \tfrac{(r-1)^{r-1}}{r^r}=\tfrac{1}{r-1}\left(1-\tfrac{1}{r}\right)^r < \tfrac{1}{e(r-1)}.$$

Finally, we consider the case $s=\omega(1)$.

\textbf{Corollary \ref{incomparable families}.\newTh{double asymptotic}\ -\ Double asymptotic version\\}
Let $k=\omega(1)$, $s=\omega(1)$ with $s \leq \frac{k}{\log_2 k}$, and $1\leq a = o(\sqrt{s})$. Then, there exist $r \sim \frac{s^a}{a!}$ pairwise incomparable families in $\mathcal{P}([k])$ of size $b \sim \frac{a^a}{e^a a!} \frac{2^k}{r}$.

The proportion of the union of all families in $\mathcal{P}([k])$ is $\frac{a^a}{e^a a!}$, which is constant for constant $a$,\linebreak and $\sim\frac{1}{\sqrt{2\pi a}}$ for $a=\omega(1)$ (Stirling's approximation).

\textbf{Proof. \ } It follows from $s \leq \frac{k}{\log_2 k} \leq k$ that 
$$\log_2 s \leq \log_2 k - \log_2 \log_2 k \leq \log_2 k - \log_2 \log_2 s,$$
and hence,
$$k \geq s \log_2 k \geq s\ (\log_2 s + \log_2 \log_2 s) = s\ (\log_2 s + \omega(1)).$$
Therefore, we can apply Theorem \ref{incomparable families}.\ref{asymptotic}. We receive $\binom{s}{a}$ pairwise incomparable families in $\mathcal{P}([k])$ of size $b\sim\frac{a^a(s-a)^{s-a}}{s^s} 2^k$.

We have to show the asymptotic equalities using $a=o(\sqrt{s})$.
$\binom{s}{a} \sim \frac{s^a}{a!}$ follows from the simple bounds $\tfrac{(s-a)^a}{a!} \leq \binom{s}{a} \leq \tfrac{s^a}{a!}$ and $(s-a)^a\sim s^a$:
$$ 1 \geq \left(\tfrac{s-a}{s}\right)^a = \left(1+\tfrac{a}{s-a}\right)^{-a} \geq e^{-a^2/(s-a)} = e^{-o(1)} \sim 1.$$ 
$\tfrac{a^a(s-a)^{s-a}}{s^s} \sim \left(\frac{a}{es}\right)^a$ follows from $\left(\tfrac{s-a}{s}\right)^{s-a} \sim e^{-a}$:

\hspace*{15mm} $e^{-a} \leq \left( 1+\tfrac{a}{s-a} \right)^{-(s-a)} = \left( \tfrac{s-a}{s }\right)^{s-a} = \left( 1-\tfrac{a}{s} \right)^{s-a} \leq e^{-a}e^{a^2/s} \sim e^{-a}.$ \hfill $\Box$

\section{Exceeding sequences of families of sets} \label{exceeding families}

In this section, we examine the maximum size $b$ of equal-sized families in arbitrarily long $d$-exceeding sequences in $\mathcal{P}([k])$. For $d=1$, we are able to determine the exact value, which is $b=\tfrac{1}{4}2^k$.
As a consequence of the following lemma, lower bounds on the size of $d+1$ pairwise incomparable families can be transferred to arbitrarily long $d$-exceeding sequences.

\textbf{Lemma \ref{exceeding families}.\newTh{transition lemma}\ -\ Incomparable and exceeding families\\}
Assume that there are $r$ pairwise incomparable families of size $b$ in $\mathcal{P}([k])$. 
Then, for all $l\in \mathbb{N}$, there is an $(r-1)$-exceeding sequence $\mathcal{B}_1,\dots,\mathcal{B}_l \ \subseteq \mathcal{P}([k])$ of families of size $b$.

\textbf{Proof. \ } Let $\mathcal{A}_1,\dots,\mathcal{A}_r \ \subseteq \mathcal{P}([k])$ be pairwise incomparable. Then, $\mathcal{B}_1,\dots,\mathcal{B}_l$ with $\mathcal{B}_i = \mathcal{A}_j$ if and only if $i= j \mod r$ is $(r-1)$-exceeding. Two numbers $i,i' \in [l]$ with $i'-i \in [r-1]$ correspond to $j\neq j'$ in $[r]$. Hence, the families $\mathcal{B}_i = \mathcal{A}_j \neq \mathcal{A}_{j'} = \mathcal{B}_{i'}$ are incomparable, in particular $\mathcal{B}_i \vdash \mathcal{B}_{i'}$. \hfill $\Box$\\

The author conjectures that Corollary \ref{incomparable families}.\ref{constant number} is best possible (which can also be derived from the conjecture in \cite{New}), even in combination with Lemma \ref{exceeding families}.\ref{transition lemma}:   

\textbf{Conjecture \ref{exceeding families}.\newTh{Conjecture}\ -\ The maximum size of \boldmath $d$-exceeding families\\}
For $l\in\mathbb{N}$ large enough, there is no $(r-1)$-exceeding sequence $\mathcal{B}_1,\dots,\mathcal{B}_l \ \subseteq \mathcal{P}([k])$ of families of size $b>\frac{(r-1)^{r-1}}{r^r}2^k$. As a consequence, there do not exist $r$ pairwise incomparable families of size $b$ in $\mathcal{P}([k])$.

Clearly, for $r=1$, Conjecture \ref{exceeding families}.\ref{Conjecture} is true. We are going to prove it for $r=2$. For this case, the second part has already been shown by Seymour in \cite{Seymour}. 
Similarly to Seymour, we use a version of Kleitman's Lemma \cite{Kleitman}. For any family $\mathcal{X} \subseteq \mathcal{P}([k])$, we define the two families $\hat{\mathcal{X}} = \{Y \subseteq [k]\ |\ \exists X \in \mathcal{X}: X \supseteq Y\}$ and $\check{\mathcal{X}} = \{Y \subseteq [k]\ |\ \exists X \in \mathcal{X}: X \subseteq Y\}$. Clearly, $\mathcal{X} \subseteq \hat{\mathcal{X}} \cap \check{\mathcal{X}}$.

\textbf{Lemma \ref{exceeding families}.\newTh{Kleitman-like}\ -\ A version of Kleitman's Lemma  \cite{Kleitman}\\ }
For every $\mathcal{X} \subseteq \mathcal{P}([k])$, \vspace*{-2mm}
$$ \frac{|\hat{\mathcal{X}}|}{2^k} \ \cdot \ \frac{|\check{\mathcal{X}}|}{2^k}\ \geq \ \frac{|\hat{\mathcal{X}}\cap \check{\mathcal{X}}|}{2^k}\ \geq \ \frac{|\mathcal{X}|}{2^k} \ .$$
\begin{figure}[H]
\begin{center}
	\vspace*{-4mm}
	\includegraphics[width=28mm]{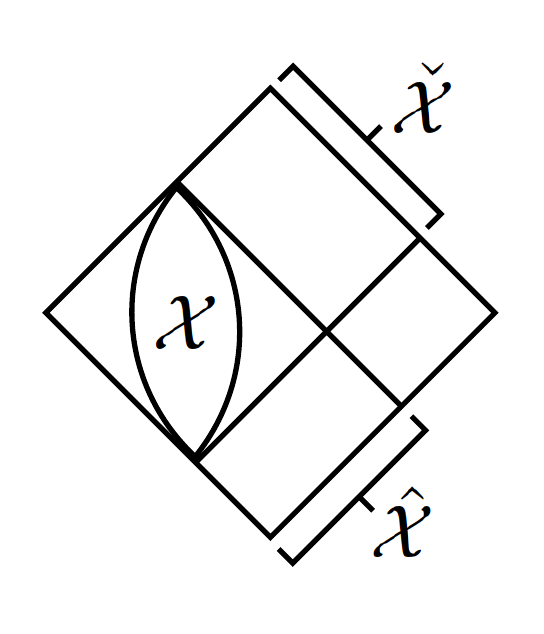}
	\vspace*{-4mm}
	\caption{As an analogy, we can think of the three fractions as areas of rectangles in the union square, which corresponds to $\mathcal{P}([k])$. The families $\hat{\mathcal{X}}$ and $\check{\mathcal{X}}$ correspond to rectangles with one side of length $1$ and one of length equal to their area $\frac{|\hat{\mathcal{X}}|}{2^k}$ or $\frac{|\check{\mathcal{X}}|}{2^k}$, respectively. Hence, the area of their intersection is equal to $\frac{|\hat{\mathcal{X}}|}{2^k} \cdot \frac{|\check{\mathcal{X}}|}{2^k}$.}
	\label{picture lemma}
\end{center}
\end{figure} \vspace*{-7mm}

To use Lemma \ref{exceeding families}.\ref{Kleitman-like} on exceeding sequences of families, we also need some relation between consecutive families in the sequence. This is provided by the following lemma.

\textbf{Lemma \ref{exceeding families}.\newTh{exceeding Lemma}\ -\ A relation between consecutive families in exceeding sequences\\ }
Let $\mathcal{X} \vdash \mathcal{Y}$ be families in $\mathcal{P}([k])$. Then, $\hat{\mathcal{X}}$ and $\check{\mathcal{Y}}$ are disjoint, and hence,
$$\frac{|\hat{\mathcal{X}}|}{2^k} + \frac{|\check{\mathcal{Y}}|}{2^k}\ \leq\ 1.$$
\textbf{Proof. \ } Assume that there is a set $C$ in both families. $C\in \hat{\mathcal{X}}$ implies that $C \subseteq X$ for some $X \in \mathcal{X}$. $C\in \check{\mathcal{Y}}$ implies that $C \supseteq Y$ for some $Y \in \mathcal{Y}$. This yields $X \supseteq C \supseteq Y$, which contradicts $\mathcal{X} \vdash \mathcal{Y}$. \hfill $\Box$\\ 

The combination of both lemmas can be used to prove the non-existence of various exceeding sequences of  families of given sizes. We only want to consider equal sizes. The following theorem should serve as a warm up:

\textbf{Theorem \ref{exceeding families}.\newTh{l4}\ -\ Exceeding families of size \boldmath $\big\lceil\tfrac{1}{3}2^k\big\rceil$\\}
There is no exceeding sequence $\mathcal{B}_1\vdash\mathcal{B}_2\vdash\mathcal{B}_3\vdash\mathcal{B}_4 \ \subseteq \mathcal{P}([k])$ of families of size $\big\lceil\tfrac{1}{3}2^k\big\rceil$.

\textbf{Proof. \ } Assume that $\mathcal{B}_1\vdash\mathcal{B}_2\vdash\mathcal{B}_3\vdash\mathcal{B}_4 \ \subseteq \mathcal{P}([k])$ of families of size $\big\lceil\tfrac{1}{3}2^k\big\rceil$ is exceeding. Note that $\big\lceil\tfrac{1}{3}2^k\big\rceil > \tfrac{1}{3}2^k$. We conclude from the previous lemmas that 
$$\frac{|\hat{\mathcal{B}}_2|}{2^k} \geq \frac{|\mathcal{B}_2|}{2^k} \left( \frac{|\check{\mathcal{B}}_2|}{2^k} \right)^{-1} > \frac{1}{3} \left( 1 - \frac{|\hat{\mathcal{B}}_1|}{2^k} \right)^{-1} \geq \frac{1}{3} \left( 1 - \frac{|\mathcal{B}_1|}{2^k} \right)^{-1} > \frac{1}{3} \left( 1 - \frac{1}{3}\right)^{-1} = \frac{1}{2},$$
and similarly, \vspace*{-2mm}
$$\frac{|\check{\mathcal{B}}_3|}{2^k} \geq \frac{|\mathcal{B}_3|}{2^k} \left( \frac{|\hat{\mathcal{B}}_3|}{2^k} \right)^{-1} > \frac{1}{3} \left( 1 - \frac{|\check{\mathcal{B}}_4|}{2^k} \right)^{-1} \geq \frac{1}{3} \left( 1 - \frac{|\mathcal{B}_4|}{2^k} \right)^{-1} > \frac{1}{3} \left( 1 - \frac{1}{3}\right)^{-1} = \frac{1}{2}.$$
But this contradicts $\frac{|\hat{\mathcal{B}}_2|}{2^k} + \frac{|\check{\mathcal{B}}_3|}{2^k} \leq 1$. \hfill $\Box$\\ 

We now consider exceeding sequences of families of size $\tfrac{1}{4}2^k$, which exist by Lemma \ref{exceeding families}.\ref{transition lemma}. But if such a sequence is sufficiently long, one cannot "squeeze" an additional set into the center of the sequence:

\textbf{Theorem \ref{exceeding families}.\newTh{r2}\ -\ Exceeding families of size \boldmath $\tfrac{1}{4}2^k$\\}
Let $k\geq 2$ and $s=\lfloor 2^{(k-1)/2} \rfloor$. For every exceeding sequence $\mathcal{B}_1\vdash\dots\vdash\mathcal{B}_{4s} \ \subseteq \mathcal{P}([k])$ of families of size $\frac{1}{4}2^k$, there is no $B\subseteq [k]$ such that $\mathcal{B}_{2s} \vdash \{B\} \vdash \mathcal{B}_{2s+1}$.\\
The length of this sequence is $4s \leq 4\cdot 2^{(k-1)/2} = 2^{(k+3)/2}$.

\textbf{Proof. \ } We prove inductively that \vspace*{-2mm}
\begin{equation}
\frac{|\hat{\mathcal{B}}_i|}{2^k} \geq \frac{1}{2} - \frac{1}{2(i+1)}
\end{equation}
For $i=1$, this is clearly true since $\frac{|\hat{\mathcal{B}}_1|}{2^k} \geq \frac{|\mathcal{B}_1|}{2^k} = \frac{1}{4} =\frac{1}{2} - \frac{1}{2(1+1)}$.
For $i>1$, we conclude using the previous lemmas (and the induction hypothesis for $i-1$)
$$\frac{|\hat{\mathcal{B}}_i|}{2^k} \geq \frac{|\mathcal{B}_i|}{2^k} \left( \frac{|\check{\mathcal{B}}_i|}{2^k} \right)^{-1} \geq \frac{1}{4} \left(1-\frac{|\hat{\mathcal{B}}_{i-1}|}{2^k} \right)^{-1} \geq \frac{1}{4} \left(1 - \frac{1}{2} +\frac{1}{2i} \right)^{-1}$$ \vspace*{-2mm}
$$ = \frac{1}{4} \left(\frac{i+1}{2i} \right)^{-1} = \frac{i}{2(i+1)} = \frac{1}{2} - \frac{1}{2(i+1)}.$$
This inequality works well for small $i$. For large $i$, we use instead that 
\begin{equation}
|\hat{\mathcal{B}}_{i+1}|>|\hat{\mathcal{B}}_i|\,, \text{\quad unless } |\hat{\mathcal{B}}_i| =\tfrac{1}{2}2^k\, \text{,\quad since}
\end{equation}
$$\frac{|\hat{\mathcal{B}}_i|}{|\hat{\mathcal{B}}_{i+1}|} = \frac{|\hat{\mathcal{B}}_i|}{2^k} \frac{2^k}{|\hat{\mathcal{B}}_{i+1}|} \leq \frac{|\hat{\mathcal{B}}_i|}{2^k} \frac{2^k}{|\mathcal{B}_{i+1}|} \frac{|\check{\mathcal{B}}_{i+1}|}{2^k} \leq 4\ \frac{|\hat{\mathcal{B}}_i|}{2^k} \left(1- \frac{|\hat{\mathcal{B}}_i|}{2^k}\right)  = 1-4 \left(\frac{1}{2}-\frac{|\hat{\mathcal{B}}_i|}{2^k}\right)^2 < 1.$$
Note that the sizes are natural numbers, so their difference is at least 1.

Now, we can conclude that $|\hat{\mathcal{B}}_{2s}| \geq \frac{1}{2} 2^k$. By (1), we have that 
$$\frac{|\hat{\mathcal{B}}_s|}{2^k} \geq \frac{1}{2} - \frac{1}{2(s+1)} > \frac{1}{2} - \frac{1}{2\cdot 2^{(k-1)/2}}.$$
We obtain $|\hat{\mathcal{B}}_s| > \frac{1}{2} 2^k - 2^{(k-1)/2}$, so even $|\hat{\mathcal{B}}_s| \geq \lceil \frac{1}{2} 2^k - 2^{(k-1)/2} \rceil = \frac{1}{2} 2^k - s$,\, so by (2), $|\hat{\mathcal{B}}_{2s}| \geq |\hat{\mathcal{B}}_s| + s \geq  \frac{1}{2} 2^k$. Similarly, one can show that $|\check{\mathcal{B}}_{2s+1}| \geq \frac{1}{2} 2^k$.

By Lemma \ref{exceeding families}.\ref{exceeding Lemma}, the families $\hat{\mathcal{B}}_{2s}$ and $\check{\mathcal{B}}_{2s+1}$ are disjoint, so they partition $\mathcal{P}([k])$. Therefore, every $B\subseteq [k]$ is contained in one of the two families, and this contradicts $B_{2s} \not\supseteq B$ or $B \not\supseteq B_{2s+1}$ for some $B_{2s}\in\mathcal{B}_{2s}$ or $B_{2s+1}\in\mathcal{B}_{2s+1}$. \hfill $\Box$\\

By attaching the conflicting additional set $B\notin \mathcal{B}_{2s}$ to $\mathcal{B}_{2s}$, we conclude the following:

\textbf{Corollary \ref{exceeding families}.\newTh{r2 corollary}\ -\ The maximum size of exceeding families \\}
In $\mathcal{P}([k])$, there are no arbitrary long exceeding sequences of families where  one family in the middle of the sequence has size $\frac{1}{4}2^k+1$ and all others have size $\frac{1}{4}2^k$.

Hence, in $\mathcal{P}([k])$, there are no two incomparable families of the sizes $\frac{1}{4}2^k+1$ and $\frac{1}{4}2^k$, which agrees with Seymour's result in \cite{Seymour}.

We finally prove a result on $d$-exceeding sequences for a general $d$. Unfortunately, the bound we obtain is worse than the one from Conjecture \ref{exceeding families}.\ref{Conjecture}. In \cite{On a conjecture of Hilton} and in \cite{New}, this is done for $r$ pairwise incomparable families (of not necessarily distinct sizes). The general idea is to unite families to trace back the analysis to the case $r=2$.

\textbf{Theorem \ref{exceeding families}.\newTh{general r}\ -\ A bound on the size of \boldmath $d$-exceeding families \vspace*{0.5ex}\\}
Let~~ $\beta = \left\{ \begin{aligned} \frac{1}{2r}& \ , &\quad r \text{ even\,,} \\ r-\sqrt{r^2-1}\ =\ \frac{1}{2r}&\left(1+\frac{1+o(1)}{4r^2} \right) , &\quad r \text{ odd\,.} \\ \end{aligned} \right.$

Then, for an $l$ large enough, there is no $(r-1)$-exceeding sequence of length $l$ of families of size $>\beta 2^k$ in $\mathcal{P}([k])$. Hence, there are no $r$ pairwise incomparable families of size $>\beta 2^k$ in $\mathcal{P}([k])$.

\textbf{Proof.\ } Let $\mathcal{B}_1,\dots,\mathcal{B}_l \ \subseteq \mathcal{P}([k])$ of size $>\beta 2^k$ be $(r-1)$-exceeding.

(i)\ Let $r=2s$. Unite $s$ consecutive families to larger families $\mathcal{A}_1,\dots,\mathcal{A}_{l'} \ \subseteq \mathcal{P}([k])$ with $l'=\lfloor l/s \rfloor$ and $\mathcal{A}_i := \mathcal{B}_{(i-1)s+1} \cup \dots \cup \mathcal{B}_{is}$. Each $\mathcal{A}_i$ has size $>s\beta 2^k=\tfrac{r\beta}{2}2^k=\tfrac{1}{4}2^k$ and their sequence is exceeding. Using Lemma \ref{exceeding families}.\ref{Kleitman-like} and \ref{exceeding families}.\ref{exceeding Lemma}, we prove that $|\hat{\mathcal{A}}_1|<\dots<|\hat{\mathcal{A}}_{l'}| \leq 2^k$, which yields a contradiction for $l'>2^k$. Indeed, 
$$\frac{|\hat{\mathcal{A}}_i|}{|\hat{\mathcal{A}}_{i+1}|} = \frac{|\hat{\mathcal{A}}_i|}{2^k} \frac{2^k}{|\hat{\mathcal{A}}_{i+1}|} \leq \frac{|\hat{\mathcal{A}}_i|}{2^k} \frac{2^k}{|\mathcal{A}_{i+1}|} \frac{|\check{\mathcal{A}}_{i+1}|}{2^k} < 4\ \frac{|\hat{\mathcal{A}}_i|}{2^k} \left(1- \frac{|\hat{\mathcal{A}}_i|}{2^k}\right)  = 1-4 \left(\frac{1}{2}-\frac{|\hat{\mathcal{A}}_i|}{2^k}\right)^2 \leq 1.$$ \newpage

(ii)\ Let $r=2s+1$. Unite $s$ consecutive families to larger families, $\mathcal{A}_1,\dots,\mathcal{A}_{l'} \ \subseteq \mathcal{P}([k])$ with $l'=\lfloor l/(s+1) \rfloor$ and
$\mathcal{A}_i := \mathcal{B}_{(i-1)(s+1)+1} \cup \dots \cup \mathcal{B}_{i(s+1)-1}$. Each $\mathcal{A}_i$ has size $> s\beta 2^k$. 

We have that $\mathcal{A}_i \vdash \mathcal{B}_{i(s+1)} \vdash \mathcal{A}_{i+1}$ and $\mathcal{A}_i \vdash \mathcal{A}_{i+1}$. Hence, by Lemma \ref{exceeding families}.\ref{exceeding Lemma}, the three families $\hat{\mathcal{A}}_i$, $\check{\mathcal{A}}_{i+1}$ and $\mathcal{B}_{i(s+1)} \subseteq \check{\mathcal{B}}_{i(s+1)}\cap \hat{\mathcal{B}}_{i(s+1)}$ are pairwise disjoint, and hence
$$ 1\ \geq\ \frac{|\hat{\mathcal{A}}_i|}{2^k} + \frac{|\check{\mathcal{A}}_{i+1}|}{2^k} + \frac{|\mathcal{B}_{i(s+1)}|}{2^k}\ >\ \frac{|\hat{\mathcal{A}}_i|}{2^k} + \frac{|\check{\mathcal{A}}_{i+1}|}{2^k} + \beta\ .$$
Again, we prove that $|\hat{\mathcal{A}}_1|<\dots<|\hat{\mathcal{A}}_{l'}| \leq 2^k$, which yields a contradiction for $l'>2^k$.
$$\frac{|\hat{\mathcal{A}}_i|}{|\hat{\mathcal{A}}_{i+1}|} = \frac{|\hat{\mathcal{A}}_i|}{2^k} \frac{2^k}{|\hat{\mathcal{A}}_{i+1}|} \leq \frac{|\hat{\mathcal{A}}_i|}{2^k} \frac{2^k}{|\mathcal{A}_{i+1}|} \frac{|\check{\mathcal{A}}_{i+1}|}{2^k} < \frac{1}{s\beta}\ \frac{|\hat{\mathcal{A}}_i|}{2^k} \left(1- \beta - \frac{|\hat{\mathcal{A}}_i|}{2^k}\right)$$ \vspace*{-2mm}
$$  = \frac{1}{s\beta} \left(\frac{(1-\beta)^2}{4}-\left( \frac{1-\beta}{2} - \frac{|\hat{\mathcal{A}}_i|}{2^k} \right)^2 \right) \leq \frac{(1-\beta)^2}{4s\beta} = 1.$$
To satisfy the last equality, we have chosen $\beta=r-\sqrt{r^2-1}$ to be a solution of\vspace*{1mm} \\ $0=\beta^2-2r\beta+1=\beta^2-(4s+2)\beta+1$. 

Finally, we have to prove that indeed $\beta = r-\sqrt{r^2-1}\ =\ \frac{1}{2r}\left(1+\frac{1+o(1)}{4r^2} \right)$. We have that $\beta = \tfrac{1}{2r}(1+\beta^2)$. From this and $\beta\leq 1$, we conclude one after another that $\beta = o(1)$, hence $\beta = \tfrac{1}{2r}(1+o(1))$, and hence $\beta = \frac{1}{2r}\left(1+\frac{1+o(1)}{4r^2} \right)$. \hfill $\Box$

\section{Exceeding sequences of sets} \label{exceeding sets}

In this section, we consider $d$-exceeding sequences of sets, which can be looked upon as families of size $1$. The main goal is to estimate $\delta(k)$, which we have defined as the maximum $d$ such that there exist arbitrarily long $d$-exceeding sequences of subsets of $[k]$. The following lemma is used to derive some lower bounds on $\delta(k)$ from results of the previous sections.

\textbf{Lemma \ref{exceeding sets}.\newTh{exceeding families and subsets}\ -\ \boldmath Exceeding sets from exceeding families\\}
If there is a $d$-exceeding sequence $\mathcal{B}_1,\dots,\mathcal{B}_l \ \subseteq \mathcal{P}([k])$ of families of size $b$, then there is a $db$-exceeding sequence $B_1,\dots,B_{lb} \subseteq [k]$.

\textbf{Proof.\ } Choose $B_1,\dots,B_{lb} \subseteq [k]$ such that, for every $i\in[l]$, $\mathcal{B}_i=\{B_{(i-1)b+1},\dots,B_{ib}\}$ is ordered by non-decreasing size, that is, $|B_{(i-1)b+1}|\leq |B_{(i-1)b+2}|\leq \dots\leq |B_{ib}|$.

Now let $j,j'\in [lb]$ with $j'-j\in [db]$. There exist $i,i'\in [l]$ and $a,a'\in [b]$ such that $j=(i-1)b+a$ and $j'=(i'-1)b+a'$. It follows that $0\leq i'-i \leq d$.

If $i'-i \in [d]$, we obtain $B_j\not\supseteq B_{j'}$ from $\mathcal{B}_i\vdash \mathcal{B}_{i'}$. Otherwise $i=i'$ and $a<a'$, so $B_j\not\supseteq B_{j'}$ follows from $B_j\neq B_{j'}$ and $|B_j|\leq |B_{j'}|$. \hfill $\Box$\\

\textbf{Theorem \ref{exceeding sets}.\newTh{Lower bounds}\ -\ Lower bounds on \boldmath $\delta(k)$}
\begin{enumerate}
	\item[(i)] $\delta(k) \geq \binom{k}{\lfloor k/2 \rfloor} -1$. \vspace*{-2mm}
	\item[(ii)] For $k\geq 4$, \ $\delta(k) \geq \frac{5}{16} 2^k$. \vspace*{-2mm}
	\item[(iii)] If $r$ is a power of $2$ and $k\geq r \log_2 r$, then $\delta(k) \geq \left( 1-\tfrac{1}{r} \right)^r 2^k \geq \left( 1-\tfrac{1}{r} \right) \tfrac{1}{e} 2^k$. \vspace*{-2mm}
	\item[(iv)] $\delta(k) \geq \left( 1-\tfrac{2 \log_2 k}{k} \right) \tfrac{1}{e} 2^k \sim \tfrac{1}{e} 2^k$.
\end{enumerate}
\textbf{Proof. \ } We use Lemma \ref{exceeding sets}.\ref{exceeding families and subsets} combined with Lemma \ref{exceeding families}.\ref{transition lemma}, that is, if there are $r$ pairwise incomparable families of size $b$ in $\mathcal{P}([k])$, then $\delta(k)\geq (r-1)b$.
\begin{enumerate}
	\item[(i)] The $\binom{k}{\lfloor k/2 \rfloor}$ subsets of size $\lfloor k/2 \rfloor$ of $[k]$ are pairwise incomparable. (Use them as families of size $1$.)
	\item[(ii)] Use the six pairwise incomparable families
$$\{B\in [k]:\ B\cap \{1,2,3,4\} = A\}, \text{\quad for } A\subseteq \{1,2,3,4\} \text{ of size } 2.$$
Each of the families has size $\tfrac{1}{16} 2^k$. The construction of the six families corresponds to Theorem \ref{incomparable families}.\ref{exact} with $s=4$ and $a=2$.
	\item[(iii)] By Theorem \ref{incomparable families}.\ref{exact} with $s=r$ and $a=1$, there exist $r$ pairwise incomparable families of size $b=\tfrac{(r-1)^{r-1}}{r^r} 2^k$. It holds that $(r-1)b = \left( \tfrac{r-1}{r} \right)^r 2^k = \left( 1-\tfrac{1}{r} \right)^r 2^k$. We can bound further
	$\left( \tfrac{r-1}{r} \right)^{r-1} = \left( 1+\tfrac{1}{r-1} \right)^{-(r-1)} \geq e^{-\tfrac{r-1}{r-1}} = \tfrac{1}{e}$.
	\item[(iv)] In (iii), we may take the maximal number $r$ that is a power of 2 under the condition that $k\geq r\log_2 r$. We have to show that $r\geq \tfrac{k}{2 \log_2 k}$, and then the claim follows from (iii). Note that this is true for $k\leq 7$, where $r=2$. For $k\geq 8$, we have $k\geq r \log_2 r \geq 2r$. By the maximality of $r$, we obtain $k < 2r \log_2(2r) \leq 2r \log_2 k$.   \hfill $\Box$
\end{enumerate} 

(i) yields $\delta(2)\geq 1$, $\delta(3)\geq 2$ and $\delta(4)\geq 5$. In the following theorems, we are going to prove that those bounds are attained. For $k\geq 5$, (ii) yields the better lower bounds $\delta(5)\geq 10$, $\delta(6)\geq 20$, $\delta(7)\geq 40$. For $k\geq 8$, (iii) yields even better bounds, starting with $\delta(8)\geq 81$.

Conjecture \ref{exceeding families}.\ref{Conjecture} implies that, for $1=b>\frac{d^d}{(d+1)^{d+1}}2^k$ and $l\in\mathbb{N}$ large enough, there is no $d$-exceeding sequence $B_1,\dots,B_l \ \subseteq [k]$:

\textbf{Conjecture \ref{exceeding sets}.\newTh{Conjecture2}\ -\ \boldmath $\delta(k)\sim\tfrac{1}{e}2^k$}
$$\delta(k) \leq \max\left\{d\in \mathbb{N}: 1\leq \frac{d^d}{(d+1)^{d+1}}2^k\right\} 
< \frac{1}{e}2^k\ .$$
For the right inequality, note that, if $d$ satisfies $1\leq \frac{d^d}{(d+1)^{d+1}}2^k$, then indeed 
$$d \leq \frac{d^{d+1}}{(d+1)^{d+1}} = \left(1-\frac{1}{d+1}\right)^{d+1} < \frac{1}{e}\,.$$

In the remaining section, we want to present some upper bounds on $\delta(k)$. The first one follows from considering a sequence of subsets of $[k]$ of length $2^k+1$, which is the lowest interesting length. A sequence of length $2^k$ which is $d$-exceeding for all $d\in\mathbb{N}$ is obtained by using all $2^k$ distinct subsets of $[k]$ ordered by non-decreasing size. 

\textbf{Theorem \ref{exceeding sets}.\newTh{large d}\ -\ Sequences of length \boldmath $2^k+1$\\} 
There exists a $d$-exceeding sequence of subsets of $[k]$ of length $2^k+1$ if and only if\\ $d \leq 2^k -2^{\lfloor k/2 \rfloor} - 2^{\lceil k/2 \rceil} +1$.

\textbf{Proof. \ } Let $d = 2^k -2^{\lfloor k/2 \rfloor} - 2^{\lceil k/2 \rceil} +1$. For some fixed set $A\subseteq [k]$ of size $\lfloor k/2 \rfloor$, consider the families 
$$\hat{A} = \{B \subseteq [k]: B\subseteq A\}, \quad \check{A} = \{B \subseteq [k]: B\supseteq A\}\ \text{ and }\ \bar{A} = \mathcal{P}([k]) \setminus (\hat{A} \cup \check{A})$$
of the sizes $2^{\lfloor k/2 \rfloor}$,\ $2^{\lceil k/2 \rceil}$\ and\ $2^k -2^{\lfloor k/2 \rfloor} - 2^{\lceil k/2 \rceil} +1 =d$, respectively. Note that both the families $\hat{A}$ and $\check{A}$ contain the set $A$.

We have that $\hat{A} \vdash \bar{A} \vdash \check{A}$: From $B \subseteq A$ and $B' \not\subseteq A$, it follows that $B \not\supseteq B'$. Similarly, from $B \not\supseteq A$ and $B' \supseteq B$, it follows that $B \not\supseteq B'$.
A $d$-exceeding sequence of subsets of $[k]$ is obtained like in the proof of Lemma \ref{exceeding sets}.\ref{exceeding families and subsets}.

Consider an arbitrary sequence $B_1,\dots,B_{2^k+1} \subseteq [k]$. Since two sets must coincide, there exists some $d\geq 0$ and some $s \in [2^k-d]$ with $B_{s} \supseteq B_{s+d+1}$. We choose the minimal such $d$, that is, the sequence is $d$-exceeding, but not $(d+1)$-exceeding. We have to show $d \leq 2^k -2^{\lfloor k/2 \rfloor} - 2^{\lceil k/2 \rceil} +1$.

Let $\mathcal{B}:=\{B_{s+1},B_{s+2},\dots,B_{s+d}\}$. $\mathcal{B}$ contains $d$ pairwise distinct subsets of $[k]$ and for each $B \in \mathcal{B}$, $B_{s} \not\supseteq B \not\supseteq B_{s+d+1} $. Since $B_{s} \supseteq B_{s+d+1}$, all $B \in \mathcal{B}$ are incomparable to $B_{s}$. This implies
$d = |\mathcal{B}| \leq 2^k - 2^{|B_{s}|} - 2^{k-|B_{s}|} + 1 \leq 2^k -2^{\lfloor k/2 \rfloor} - 2^{\lceil k/2 \rceil} +1$. We have used that, for $0\leq s \leq t \leq k$ with $s+t=k$, it holds that $2^s + 2^t \geq 2^{\lfloor k/2 \rfloor} + 2^{\lceil k/2 \rceil}$: If $t = \lceil k/2 \rceil$, we have equality. If $t \geq \lceil k/2 \rceil +1$, we obtain $2^s + 2^t > 2^t \geq 2^{\lceil k/2 \rceil +1} \geq 2^{\lfloor k/2 \rfloor} + 2^{\lceil k/2 \rceil}$. \hfill $\Box$\\

In particular, we obtain that there is no $2$-exceeding sequence of subsets of $[2]$ of length $5$, and hence $\delta(2)=1$. Moreover, there is no $4$-exceeding sequence of subsets of $[3]$ of length $9$, and hence $\delta(3)\leq 3$. In the following two theorems, we prove that $\delta(3)=2$ and $\delta(4)=5$.

\textbf{Theorem \ref{exceeding sets}.\newTh{d3}\ -\ \boldmath $\delta(3)=2$.\\}
There is a $3$-exceeding sequence of subsets of $\{1,2,3\}$ of length $l$ if and only if $l\leq 10$.

\textbf{Proof. \ } The following sequence of length $10$ is $3$-exceeding:
$$\emptyset, \{1\}, \{2\}, \{3\}, \{2,3\}, \{1\}, \{1,2\}, \{1,3\}, \{2,3\}, \{1,2,3\}.$$
Let $B_1,\dots,B_l \subseteq \{1,2,3\}$ be $3$-exceeding. We have to prove that $l\leq 10$.
Each of the sets $B_2,\dots,B_{l-1}$ contains one or two elements. Let $h$ be the greatest number such that $|B_h|= 1$. We have that $h\geq l-4$, since otherwise the sets $B_{l-4},B_{l-3},B_{l-2},B_{l-1}$ each contain two elements and are distinct. But there are only three such sets. \enlargethispage{5mm}

Assume that $h=7$. We may assume that $B_7=\{1\}$. Then, $B_4,B_5,B_6$ are non-empty and do not contain $1$ as an element. The only such sets are $\{2\}, \{3\}, \{2,3\}$. Hence, $B_4,B_5$ are the two sets $\{2\}, \{3\}$. It follows that $B_2\neq B_3$ are non-empty and do not contain $2$ or $3$ as elements, but $\{1\}$ is the only such set. This contradictions occurs also if $h>7$. So $h\leq 6$, and $l \leq h+4 \leq 10$. \hfill $\Box$\\

\textbf{Theorem \ref{exceeding sets}.\newTh{d4}\ -\ \boldmath $\delta(4)=5$.\\}
There is a $6$-exceeding sequence of subsets of $\{1,2,3,4\}$ of length $l$ if and only if $l\leq 24$.

\textbf{Proof. \ } The following sequence of length $24$ is $4$-exceeding (non-empty sets are written as words, for example, $12 = \{1,2\}$):
$$ \emptyset, 1, 2, 3, 4, 23, 24, 34, 1, 12, 13, 14, 23, 24, 34, 234, 12, 13, 14, 123, 124, 134, 234, 1234.$$
Let $B_1,\dots,B_l \subseteq \{1,2,3,4\}$ be $6$-exceeding. We have to prove that $l\leq 24$. Each of the sets $B_2,\dots,B_{l-1}$ contains one, two or three elements. Let $h$ be the greatest index with $|B_h|= 1$ and let $h'$ be the lowest index with $|B_{h'}|= 3$.

First, consider the case $h<h'$. It follows that $h'\leq h+7$ since otherwise the seven distinct sets $B_{h+1},\dots,B_{h+7}$ must contain exactly two elements, but there are only six of those. 

Assume that $h=10$. We may assume that $B_{10}=\{1\}$. Since  $h<h'$, the six sets $B_4, \dots, B_9$ contain one or two elements and do not contain $1$ as an element. So they are exactly the six sets $\{2\}, \{3\}, \{4\}, \{2,3\}, \{2,4\}, \{3,4\}$. The three sets $\{2\}, \{3\}, \{4\}$ must be contained among $B_4, B_5, B_6, B_7$. Hence, $B_2 \neq B_3$ are non-empty and do not contains $2$, $3$ or $4$ as element, but $\{1\}$ is the only such set. This contradiction also occurs if $h>10$. It follows that $h\leq 9$.
Similarly, we can prove that $h'\geq l-8$. All together, $l \leq h'+8 \leq h+15 \leq 24$. 

Now, consider the case $h'\leq h-1$. Assume that $h=12$. We may assume that $B_{12}=\{1\}$. The six sets $B_6,\dots,B_{11}$ have to be chosen from $\{2\}, \{3\}, \{4\}, \{2,3\}, \{2,4\}, \{3,4\}, \{2,3,4\}$. It follows that, among $B_6$, $B_7$, $B_8$, there are two sets of size one, without loss of generality $\{3\}$ and $\{4\}$. Hence, $B_2, B_3, B_4, B_5$ are non-empty and do not contain $3$ or $4$ as elements, but $\{1\}, \{2\}, \{1,2\}$ are the only such sets. This contradiction also occurs if $h>12$. It follows that $h\leq 11$. Similarly, we can prove that $h'\geq l-10$. All together, $l \leq h'+10 \leq h+9 \leq 20$. \hfill $\Box$\\

Finally, we prove a general upper bound on $\delta(k)$. From Corollary \ref{exceeding families}.\ref{r2 corollary}, we could directly deduce $\delta(k)\leq \tfrac{1}{2}2^k-1$. Instead, we use the original Theorem \ref{exceeding families}.\ref{r2} to improve this bound by $1$.\linebreak Nevertheless, the improved bound is probably only attained for $k=3$ (at least not for $k=4$). 

\textbf{Theorem \ref{exceeding sets}.\newTh{upper bound}\ -\ \boldmath $\delta(k)\leq \tfrac{1}{2}2^k-2$.\\}
Let $k\geq 3$. For $l$ large enough ($l\sim 2^{3k/2}$), there is no $(\tfrac{1}{2}2^k-1)$-exceeding sequence $B_1,\dots,B_l \subseteq [k]$.

\textbf{Proof.\ } Let $s=\lfloor 2^{(k-1)/2} \rfloor$ as in Theorem \ref{exceeding families}.\ref{r2}, and $b=\tfrac{1}{4}2^k$. Assume that there is a $(2b-1)$-exceeding sequence $B_1,\dots,B_l \subseteq [k]$, where $l=(4s+2)b$.

Consider an index $c\in \{l/2 -(b-1),\dots,l/2 +(b-1)\}$ and the sequence of families $\mathcal{B}_{-2s},\dots,\mathcal{B}_{-1},\mathcal{B}_{1},\dots,\mathcal{B}_{2s}$ of size $b$, where $\mathcal{B}_{\pm t} = \lbrace B_{c\pm( (t-1)b+i )}:\ i\in [b]\,\rbrace$.

Assume that $B_{c-b}\not\supseteq B_{c+b}$. Then, the sequence of families is exceeding and additionally $\mathcal{B}_{-1} \vdash \{B_c\} \vdash \mathcal{B}_1$. But this is not possible by Theorem \ref{exceeding families}.\ref{r2}. Hence $B_{c-b} \supseteq B_{c+b}$. This is true for every index $c\in \{l/2 -(b-1),\dots, l/2 +(b-1)\}$.

We prove that the $2b$ sets $B_{l/2},\dots,B_{l/2 +2b-1}$ are pairwise incomparable. 
Let $i<j$ in $\{l/2,\dots,l/2 +2b-1\}$. We have that $B_i \not\supseteq B_j$, since the original sequence of sets is $(2b-1)$-exceeding. For the same reason, $B_{j-2b} \not\supseteq B_i$. Using our result from above with $c=j-b$, we obtain $B_{j-2b}\supseteq B_j$. Hence, also $B_j \not\supseteq B_i$.

By Sperner's Theorem, $\binom{k}{\lfloor k/2 \rfloor} \geq 2b = \tfrac{1}{2} 2^k$, which is only true for $k\leq 2$. Hence, for $k\geq 3$, we have a contradiction. \hfill $\Box$\\


\begin{thebibliography}{}

\bibitem{Seymour}
P. D. Seymour (1973): On incomparable collections of sets. Mathematika 20, 208-209.
https://doi.org/10.1112/S0025579300004782

\bibitem{false}
D. Gerbner, N. Lemons, C. Palmer, B. Patkós, V. Szécsi (2012). Cross-Sperner families. Studia Sci. Math. Hungar. 49(1), 44–51. https://doi.org/10.1556/sscmath.2011.1185

\bibitem{New}
N. Behague, A. Kuperus, N. Morrison, A. Wright (2024). Improved bounds for cross-Sperner systems. The Electronic Journal of Combinatorics 31(2). https://doi.org/10.37236/11860

\bibitem{self}
M. Krone (2024+). Cut covers of acyclic digraphs. https://arxiv.org/abs/2410.06899 

\bibitem{Sperner}
E. Sperner (1928). Ein Satz über Untermengen einer endlichen Menge. Mathematische Zeitschrift 27(1), 544-548.

\bibitem{Kleitman}
D. J. Kleitman (1966). Families of non-disjoint subsets. Journal of Combinatorial Theory 1, 153-155.

\bibitem{On a conjecture of Hilton}
Jiping Liu, Cheng Zhao (2001). On a conjecture of Hilton. Australasian Journal of Combinatorics 24, 265-274.

\end{thebibliography}
\end{document}